\documentclass{article}
\usepackage{amsmath,amssymb}
\newtheorem{df}{Definition}[section]
\newtheorem{thm}[df]{Theorem}
\newtheorem{lem}[df]{Lemma}
\title{Uncertainty product of the spherical Gauss--Weierstrass wavelet}
\author{I. Iglewska-Nowak\footnote{West Pomeranian University of Technology, Szczecin, School of Mathematics, al. Pias\-t\'ow 17, 70--310 Szczecin, Poland}}

\begin{document}

\maketitle

\bibliographystyle{amsplain}

\begin{abstract}In this paper, asymptotic behavior of the uncertainty product of Gauss-Weierstrass wavelet is investigated. It is shown that the uncertainty product is bounded from above, a feature that distinguishes Gauss-Weierstrass wavelet among other wavelet families.
\end{abstract}

\begin{bfseries}Key Words and Phrases:\end{bfseries} uncertainty product, time--frequency localization, Gauss--Weierstrass wavelet\\
\begin{bfseries}2010 MSC:\end{bfseries} 42C40

\section{Introduction}

The uncertainty product of spherical functions, a notion introduced by Narcowich and Ward in~\cite{NW96}, is a measure for the trade-off between the spatial and frequency localization. A short history of the notion and a discussion about various nomenclatures can be found in \cite[Introduction]{IIN16USW}. Similarly as in physics (Heisenberg's uncertainty principle), the uncertainty product is bounded from below. It can be interpreted as one of indicators of 'goodness' of spherical wavelets (for theory of wavelet transforms see~\cite{IIN15WT}), because well-localized wavelets yield a well-localized wavelet transform, i.e., such that 'fuzziness' of the transform reflects mostly the 'fuzziness' of the analyzed signal. In the paper~\cite{IIN16USW} I proved that the uncertainty product of a wide class of wavelets is in limit $\rho\to0$, where $\rho$ denotes the scale parameter, bounded from above by $\mathcal O(\rho^{-\alpha})$ for some $\alpha>0$. Consequently, its finite boundedness seems to be something exceptional and wavelets having this property are quite rare.

In the case of Poisson wavelets~$g_\rho^m$, $m\in\mathbb N$, over the $n$-dimensional sphere~$\mathcal S^n$ \cite{HI07,IIN15PW} -- a wavelet family which is the most popular one for applications because of their explicit representation as well as existence of discrete frames \cite{IH10,IIN16WF} -- the uncertainty product $U(g_\rho^m)$ is indeed bounded in limit $\rho\to0$ \cite{IIN16UPW}. Moreover, $\lim_{\rho\to0}U(g_\rho^m)$ approaches the optimal value for a certain sequence of $(n,m)$. 

The present paper is devoted to another very popular wavelet family, namely Gauss-Weierstrass wavelet. It is shown that its uncertainty product is bounded in limit $\rho\to0$. I was not able to compute the value of the limit with analytical tools, but in view of \cite[Theorem~3.4]{IIN16USW}, sole boundedness of the uncertainty product is a significant advantage of the wavelet family.

The paper is organized as follows. After an introduction of the necessary notions and statements in Section~\ref{sec:sphere}, the main result of the paper, Theorem~\ref{thm:uncertainty_GW} is proven in Section~\ref{sec:uncertaintyGW}.

\section{Preliminaries}\label{sec:sphere}

Let $\mathcal S$ denote the $2$-dimensional unit sphere in the Euclidean space~$\mathbb{R}^3$ and $d\sigma$ be the Lebesgue measure of~$\mathcal S$. Integrable zonal functions over the sphere (i.e., those depending only on~$\vartheta=\left<\hat e,x\right>$, where~$\hat e$ is the north-pole of the sphere $\hat e=(1,0,0)$) have the Legendre expansion
$$
f(\cos\vartheta)=\sum_{l=0}^\infty\widehat f(l)\,P_l(\cos\vartheta)
$$
with the Legendre coefficients
$$
\widehat f(l)=a_l\int_{-1}^1 f(t)\,P_l(t)\left(1-t^2\right)dt,
$$
where~$a_l$ is a constant that depends on~$l$, and~$P_l$, $l\in\mathbb N_0$, denote the Legendre polynomials.

The variances in space and momentum domain of a $\mathcal C^2(\mathcal S)$--function~$f$ with $\int_{\mathcal S}x\,|f(x)|^2\,d\sigma(x)\ne0$ are given by \cite{nLF03}
$$
\text{var}_S(f)=\left(\frac{\int_{\mathcal S}|f(x)|^2\,d\sigma(x)}{\int_{\mathcal S}x\,|f(x)|^2\,d\sigma(x)}\right)^2-1
$$
and
$$
\text{var}_M(f)=-\frac{\int_{\mathcal S}\Delta^\ast f(x)\cdot \bar f(x)\,d\sigma(x)}{\int_{\mathcal S}|f(x)|^2\,d\sigma(x)},
$$
where $\Delta^\ast$ is the Laplace--Beltrami operator on~$\mathcal S$. The quantity
\begin{equation}\label{eq:def_Uf}
U(f)=\sqrt{\text{var}_S(f)}\cdot\sqrt{\text{var}_M(f)}
\end{equation}
is called the uncertainty product (constant) of~$f$.

The uncertainty product of zonal functions may be computed from their Legendre coefficients \cite{nLF03,IIN16MR} and according to the spherical uncertainty principle it is bounded from below by~$1$, see \cite{NW96,RV97} and \cite[formula (4.37)]{GG04a}, \cite[formula (12)]{GG04b}.

\begin{lem}\label{lem:varS_varM} Let a zonal square integrable and continuously differentiable function over $\mathcal S^2$ be given by its Legendre expansion
$$
f(\cos\vartheta)=\sum_{l=0}^\infty\frac{2l+1}{4\pi}\,c_l(f)P_l(\cos\vartheta),
$$
where $\{\sqrt{2l+1}\,c_l\}\in l^2(\mathbb N_0)$. Its variances in space and momentum domain are equal to
\begin{align}
\text{var}_S(f)&=\left(\frac{\sum_{l=0}^\infty(2l+1)\,|c_l|^2}{\sum_{l=0}^\infty(l+1)\,(\overline{c_l}\,c_{l+1}+c_l\,\overline{c_{l+1}})}\right)^2-1,\label{eq:varS}\\
\text{var}_M(f)&=\frac{\sum_{l=1}^\infty l(l+1)(2l+1)\,|c_l|^2}{\sum_{l=0}^\infty(2l+1)\,|c_l|^2},\label{eq:varM}
\end{align}
whenever the series are convergent.
\end{lem}

\begin{thm}For $f\in\mathcal L^2(\mathcal S^2)\cap\mathcal C^1(\mathcal S^2)$, $U(f)\geq1$.
\end{thm}

As a consequence of~\cite[Theorem3.4]{IIN16USW} we have the following upper estimation.

\begin{thm}Let $\{\Psi_\rho\}$ be a zonal wavelet family with
$$
\widehat\Psi_\rho(l)=\frac{2l+1}{4\pi}\left[\rho^a q_\nu(l)\right]^c\,e^{-\rho^a q_\nu(l)},
$$
where $a>0$, $c>0$, and $q_{\nu}(l)=a_\nu l^\nu+a_{\nu-1}l^{\nu-1}+\dots+a_1 l+a_0$ is a polynomial of degree~$\nu$, positive and monotonously increasing for $l\geq1$. The
uncertainty product of~$\Psi_\rho$ for $\rho\to0$ behaves like
$$
U(\Psi_\rho)\leq\mathcal O\left(\rho^{\frac{-a}{2\nu}}\right).
$$
\end{thm}

\section{The uncertainty product of Gauss-Weierstrass wavelet}\label{sec:uncertaintyGW}

Gauss--Weierstrass wavelet (with respect to the measure~$\alpha(\rho)=\frac{1}{\rho}$) is given by
$$
\Psi_\rho^G(x)=\sum_{l=0}^\infty\frac{2l+1}{4\pi}\sqrt{2\rho l(l+1)}\,e^{-\rho l(l+1)}\,P_l(\cos\vartheta),\qquad x=(\vartheta,\varphi),\,\rho\in(0,\infty),
$$
see \cite[Subsec.~10.2.3]{FGS-book}.

\begin{thm}\label{thm:uncertainty_GW} The uncertainty constant of Gauss--Weierstrass wavelet is in limit $\rho\to0$ bounded by%THEOREM
\begin{equation}\label{eq:UPsi_upperbound}
 U(\Psi_\rho^G)\leq\sqrt{2\left(1+\frac{6}{e}+\frac{16}{e^2}\right)}+o(1).
\end{equation}
\end{thm}

In order to prove theorem, we need the following lemmas, which proofs are postponed to the end of the section.

\begin{lem}\label{lem:estimations}The following estimations are valid for $\rho\to0$:%LEMMA 1
\begin{align}
\tfrac{1}{8\rho^2}\!-\!\tfrac{3\sqrt 3e^{-3/2}}{8\rho\sqrt\rho}\!+\!\mathcal O(\tfrac{1}{\rho})
   &\leq\sum_{l=1}^\infty\! l(l\!+\!\tfrac{1}{2})(l\!+\!1)\,e^{\!-\!2\rho l(l\!+\!1)}
\leq\tfrac{1}{8\rho^2}\!+\!\tfrac{3\sqrt 3e^{-3/2}}{8\rho\sqrt\rho}\!+\!\mathcal O(\tfrac{1}{\rho}),\label{eq:est_v1}\\
\tfrac{1}{8\rho^2}\!-\!\tfrac{\sqrt{2\pi}\!+\!6\sqrt 3e^{\!-\!3/2}}{16\rho\sqrt\rho}\!+\!\mathcal O(\tfrac{1}{\rho})&\leq\sum_{l=1}^\infty(l\!+\!1)^2\sqrt{l(l\!+\!2)}\,e^{\!-\!2\rho(l\!+\!1)^2}
   \leq\tfrac{1}{8\rho^2}\!+\!\tfrac{3\sqrt 3e^{\!-\!3/2}}{8\rho\sqrt\rho}\!+\!\mathcal O(\tfrac{1}{\rho}),\label{eq:est_v23}\\
\tfrac{1}{8\rho^3}\!-\!\tfrac{25\sqrt 5e^{5/2}}{32\rho^2\sqrt\rho}\!+\!\mathcal O(\tfrac{1}{\rho\sqrt\rho})
   &\leq\sum_{l=1}^\infty l^2(l\!+\!\tfrac{1}{2})(l\!+\!1)^2\,e^{\!-\!2\rho l(l\!+\!1)}
   \leq\tfrac{1}{8\rho^3}\!+\!\tfrac{25\sqrt 5e^{-5/2}}{32\rho^2\sqrt\rho}\!+\!\mathcal O(\tfrac{1}{\rho\sqrt\rho}).\label{eq:est_v4}
\end{align}
\end{lem}

\begin{lem}\label{lem:estimations_varS}The series%LEMMA 2
$$
\sum_{l=1}^\infty\left[l(l+\tfrac{1}{2})(l+1)\,e^{-2\rho l(l+1)}-(l+1)^2\sqrt{l(l+2)}\,e^{-2\rho(l+1)^2}\right]
$$
is in limit $\rho\to0$ bounded from above by
$$
\left(\frac{1}{8}+\frac{3}{4e}+\frac{2}{e^2}\right)\frac{1}{\rho}+\mathcal O\left(\frac{1}{\sqrt\rho}\right).
$$
\end{lem}

\begin{bfseries}Proof of Theorem~\ref{thm:uncertainty_GW}. \end{bfseries}Substitute%PROOF OF THEOREM
$$
c_l=\sqrt{2\rho l(l+1)}\,e^{-\rho l(l+1)}
$$
to~\eqref{eq:varS} and~\eqref{eq:varM} to obtain
\begin{align*}
\text{var}_S(\Psi_\rho^G)&=\left(\frac{4\rho\sum_{l=1}^\infty l(l+\frac{1}{2})(l+1)\,e^{-2\rho l(l+1)}}
   {4\rho\sum_{l=1}^\infty(l+1)^2\sqrt{l(l+2)}\,e^{-2\rho(l+1)^2}}\right)^2-1,\\
\text{var}_M(\Psi_\rho^G)&=\frac{4\rho\sum_{l=1}^\infty l^2(l+\frac{1}{2})(l+1)^2\,e^{-2\rho l(l+1)}}
   {4\rho\sum_{l=1}^\infty l(l+\frac{1}{2})(l+1)\,e^{-2\rho l(l+1)}}.
\end{align*}
In order to estimate~$\text{var}_S(\Psi_\rho^G)$, we write it as
$$
\text{var}_S(\Psi_\rho^G)=\frac{(A+B)(A-B)}{B^2},
$$
where
\begin{align*}
A&:=\sum_{l=1}^\infty l(l+\tfrac{1}{2})(l+1)\,e^{-2\rho l(l+1)}\qquad\text{and}\\
B&:=\sum_{l=1}^\infty(l+1)^2\sqrt{l(l+2)}\,e^{-2\rho(l+1)^2}.
\end{align*}
Using estimations~\eqref{eq:est_v1} and~\eqref{eq:est_v23}, and the fact that
$$
\frac{a+\mathcal O(\rho)}{\alpha+\mathcal O(\rho)}=\frac{a}{\alpha}+o(1),\qquad\rho\to0,
$$
we conclude that
$$
\frac{A+B}{B}=\frac{\frac{1}{4\rho^2}+\mathcal O(\tfrac{1}{\rho\sqrt\rho})}{\frac{1}{8\rho^2}+\mathcal O(\tfrac{1}{\rho\sqrt\rho})}=2+o(1).
$$
On the other hand, by Lemma~\ref{lem:estimations_varS} and~\eqref{eq:est_v23},
$$
\frac{A-B}{B}\leq\frac{\left(\frac{1}{8}+\frac{3}{4e}+\frac{2}{e^2}\right)\frac{1}{\rho}+\mathcal O(\frac{1}{\sqrt\rho})}{\frac{1}{8\rho^2}+\mathcal O(\frac{1}{\rho\sqrt\rho})}
   =\left(1+\frac{6}{e}+\frac{16}{e^2}\right)\rho+o(\rho).
$$
Consequently,
$$
\text{var}_S(\Psi_\rho^G)\leq2\left(1+\frac{6}{e}+\frac{16}{e^2}\right)\rho+o(\rho).
$$
for $\rho\to0$. Similarly,
$$
\text{var}_M(\Psi_\rho^G)=\frac{1}{\rho}+o\left(\frac{1}{\rho}\right),\qquad\rho\to0,
$$
and the assertion~\eqref{eq:UPsi_upperbound} follows by definition~\eqref{eq:def_Uf}.
\hfill$\Box$

In the proofs of Lemmas~\ref{lem:estimations} and~\ref{lem:estimations_varS} the following estimation of the rest term in the trapezoidal quadrature rule will be used.

\begin{lem}\label{lem:quadrature}%LEMMA QUADRATURE
Let~$v$ be a continuous function over~$[0,\infty)$, monotonously increasing for $t\in(0,t_0)$, monotonously decreasing for $t\in(t_0,\infty)$ and such that
\begin{equation}\label{eq:0_limit_condition}
f(0)=\lim_{t\to\infty}f(t)=0.
\end{equation}
Then,
\begin{equation}\label{eq:quadrature}
\int_0^\infty v(t)\,dt-v(t_0)\leq\sum_{l=1}^\infty v(l)\leq\int_0^\infty v(t)\,dt+v(t_0).
\end{equation}
\end{lem}

\begin{bfseries}Proof. \end{bfseries}Denote by~$[t_0]$ the greatest integer less than or equal to~$t_0$. By the monotonicity of~$v$,%PROOF OF THE QUADRATURE
$$
\sum_{l=0}^{[t_0]-1}v(l)\leq\int_0^{[t_0]}v(t)\,dt\leq\sum_{l=1}^{[t_0]}v(t)
$$
and
$$
\sum_{l=[t_0]+2}^\infty v(l)\leq\int_{[t_0]+1}^\infty v(t)\,dt\leq\sum_{l=[t_0]+1}^\infty v(t).
$$
These estimations, together with
$$
\int_{[t_0]}^{[t_0]+1}v(t)\,dt\leq v(t_0)\cdot1
$$
yield the lower estimation in~\eqref{eq:quadrature}. Further, we have
$$
\min\{v([t_0]),v([t_0]+1)\}\leq\int_{[t_0]}^{[t_0]+1}v(t)\,dt
$$
and
$$
\max\{v([t_0]),v([t_0]+1)\}\leq v(t_0).
$$
Thus,
$$
v([t_0])+v([t_0]+1)\leq\int_{[t_0]}^{[t_0]+1}v(t)\,dt+v(t_0),
$$
and the upper estimation in~\eqref{eq:quadrature} follows.\hfill$\Box$

\begin{bfseries}Proof of Lemma~\ref{lem:estimations}. \end{bfseries} For the proof of estimation~\eqref{eq:est_v1} consider the function
$$
v_1(t)=t(t+\tfrac{1}{2})(t+1)e^{-2\rho t(t+1)}.
$$
It satisfies condition~\eqref{eq:0_limit_condition}. Further, its derivative equals
$$
v_1^\prime(t)=\left[3t^2+3t+\frac{1}{2}-\rho t(t+1)(2t+1)^2\right]e^{-2\rho t(t+1)}.
$$
The expression in brackets disappears for $t=t_1$ such that
\begin{equation}\label{eq:rho_t1}
\rho=\frac{6t_1^2+6t_1+1}{2t_1(t_1+1)(2t_1+1)^2},
\end{equation}
and it is nonnegative for $t<t_1$ and it is negative for $t\geq t_1$. Thus, the function~$v_1$ is monotonously increasing for $t\in[0,t_1]$ and monotonously decreasing for $t\in[t_1,\infty)$. Moreover,
$$
\frac{6t_1^2+6t_1+1}{2t_1(t_1+1)(2t_1+1)^2}=\frac{3}{4t_1^2}+\mathcal O\left(\frac{1}{t_1^3}\right)\qquad\text{for }t_1\to\infty.
$$
Multiply both sides of
$$
\rho=\frac{3}{4t_1^2}+\mathcal O\left(\frac{1}{t_1^3}\right)
$$
by~$t_1^2$ to see that
$$
t_1=\frac{\sqrt 3}{2\sqrt\rho}+\mathcal O(1)
$$
satisfies equation~\eqref{eq:rho_t1} for $\rho\to0$. Consequently,
\begin{align}
\max_{t\in\mathbb R_+}&v_1(t)=v_1(t_1)=\left(\frac{\sqrt 3}{2\sqrt\rho}+\mathcal O(1)\right)^3e^{-2\rho\left(\frac{\sqrt 3}{2\sqrt\rho}+\mathcal O(1)\right)^2}\notag\\
&=\left[\frac{3\sqrt 3}{8\rho\sqrt\rho}+\mathcal O\left(\frac{1}{\rho}\right)\right]e^{-\frac{3}{2}+\mathcal O(\sqrt\rho)}\qquad\text{for }\rho\to0.\label{eq:max_v1}
\end{align}
Since
$$
e^{-\frac{3}{2}+\mathcal O(\sqrt\rho)}=\sum_{j=0}^\infty\frac{\left(-\frac{3}{2}+\mathcal O(\sqrt\rho)\right)^j}{j!}=e^{-\frac{3}{2}}+\mathcal O(\sqrt\rho),
$$
we obtain from~\eqref{eq:max_v1},
$$
\max_{t\in\mathbb R_+}v_1(t)=\frac{3\sqrt 3}{8\rho\sqrt\rho}\,e^{-\frac{3}{2}}+\mathcal O\left(\frac{1}{\rho}\right)\qquad\text{for }\rho\to0.
$$
On the other hand, it can be checked by derivation that
$$
\int v_1(t)\,dt=-\frac{e^{-2\rho t(t+1)}(2\rho t^2+2\rho t+1)}{8\rho^2}+C,
$$
and, consequently,
$$
\int_0^\infty v_1(t)\,dt=\frac{1}{8\rho^2}.
$$
Estimation~\eqref{eq:est_v1} follows by Lemma~\ref{lem:quadrature}.

In order to obtain~\eqref{eq:est_v23} note that
$$
v_2(t)\leq(t+1)^2\sqrt{t(t+2)}e^{-2\rho(t+1)^2}\leq v_3(t)
$$
for
\begin{align*}
v_2(t)&:=t(t+1)^2\,e^{-2\rho(t+1)^2},\\
v_3(t)&:=(t+1)^3\,e^{-2\rho(t+1)^2}.
\end{align*}
Similarly as in the previous case, we obtain
\begin{align*}
v_2^\prime(t)&=(t+1)\left[3t+1-4\rho t(t+1)^2\right]e^{-2\rho(t+1)^2},\\
v_3^\prime(t)&=(t+1)^2\left[3-4\rho(t+1)^2\right]e^{-2\rho(t+1)^2},
\end{align*}
and
$$
\max_{t\in\mathbb R_+}v_2(t)=v_2(t_2),\qquad\max_{t\in\mathbb R_+}v_3(t)=v_3(t_3)
$$
with~$t_2$ such that
$$
\rho=\frac{3t_2+1}{4t_2(t_2+1)^2}=\frac{3}{4t_2^2}+\mathcal O\left(\frac{1}{t_2^3}\right)\qquad\text{for }t_2\to\infty,
$$
i.e.,
$$
t_2=\frac{\sqrt 3}{2\sqrt\rho}+\mathcal O(1)\qquad\text{for }\rho\to0,
$$
and
$$
t_3=\frac{\sqrt 3}{2\sqrt\rho}-1.
$$
Consequently,
\begin{align*}
\max_{t\in\mathbb R_+}v_2(t)&=\frac{3\sqrt 3}{8\rho\sqrt\rho}\,e^{-\frac{3}{2}}+\mathcal O\left(\frac{1}{\rho}\right)\qquad\text{for }\rho\to0,\\
\max_{t\in\mathbb R_+}v_3(t)&=\frac{3\sqrt 3}{8\rho\sqrt\rho}\,e^{-\frac{3}{2}}.
\end{align*}
Further, it follows from
\begin{align*}
\int v_2(t)\,dt&=-\frac{2[1+2\rho t(t+1)]e^{-2\rho(t+1)^2}+\sqrt{2\pi\rho}\,\text{erf}\left[\sqrt{2\rho}\,(t+1)\right]}{16\rho^2}+C\qquad\text{and}\\
\int v_3(t)\,dt&=-\frac{[1+2\rho(t+1)^2]\,e^{-2\rho(t+1)^2}}{8\rho^2}+C.
\end{align*}
that
\begin{align*}
\int_0^\infty v_2(t)\,dt&=\frac{2e^{-2\rho}-\sqrt{2\pi\rho}\,\text{erfc}(\sqrt{2\rho})}{16\rho^2}\qquad\text{and}\\
\int_0^\infty v_3(t)\,dt&=\frac{(1+2\rho)\,e^{-2\rho}}{8\rho^2}.
\end{align*}
Estimation~\eqref{eq:est_v23} is an implication of
$$
\int_0^\infty\!\!\!\!v_2(t)\,dt-\max_{t\in\mathbb R_+}v_2(t)\leq\sum_{l=1}^\infty(l+1)^2\sqrt{l(l+2)}\,e^{-2\rho(l+1)^2}
\leq\int_0^\infty\!\!\!\!v_3(t)\,dt+\max_{t\in\mathbb R_+}v_3(t)
$$
and
$$
\text{erfc}(\sqrt{2\rho})=1-2\,\sqrt\frac{2}{\pi}\sqrt\rho+\mathcal O(\rho\sqrt\rho)\qquad\text{for }\rho\to0.
$$

Finally, the derivative
$$
v_4^\prime(t)=t(t+1)\left[5t^2+5t+1-\rho t(t+1)(2t+1)^2\right]e^{-2\rho t(t+1)}
$$
of
$$
v_4(t):=t^2(t+\tfrac{1}{2})(t+1)^2\,e^{-2\rho t(t+1)}
$$
disappears for~$t=t_4$ such that
$$
\rho=\frac{5t_4^2+5t_4+1}{t_4(t_4+1)(2t_4+1)^2}=\frac{5}{4t_4^2}+\mathcal O\left(\frac{1}{t_4^3}\right)\qquad\text{for }t_4\to\infty,
$$
i.e.,
$$
t_4=\frac{\sqrt 5}{2\sqrt\rho}+\mathcal O(1))\qquad\text{for }\rho\to0.
$$
It is nonnegative for $t<t_4$ and negative for $t>t_4$. Thus,
$$
\max_{t\in\mathbb R_+}v_4(t)=v_4(t_4)=\frac{25\sqrt 5}{32\rho^2\sqrt\rho}\,e^{-\frac{5}{2}}+\mathcal O\left(\frac{1}{\rho\sqrt\rho}\right)\qquad\text{for }\rho\to0.
$$
In order to obtain~\eqref{eq:est_v4} note that
$$
\int v_4(t)\,dt=-\frac{1+2\rho t(t+1)+2\rho^2t^2(t+1)^2}{8\rho^3}\,e^{-2\rho t(t+1)}+C
$$
yields
$$
\int_0^\infty v_4(t)\,dt=\frac{1}{8\rho^3}.
$$
\hfill$\Box$

\begin{bfseries}Proof of Lemma~\ref{lem:estimations_varS}. \end{bfseries}We shall express the difference of the series as%PROOF OF LEMMA 2
\begin{align*}
\sum_{l=1}^\infty l(l+\tfrac{1}{2})&(l+1)\,e^{-2\rho l(l+1)}-\sum_{l=1}^\infty(l+1)^2\sqrt{l(l+2)}e^{-2\rho(l+1)^2}\\
&=\sum_{l=1}^\infty\left[v_1(l)+v_2(l)+v_3(l)\right]
\end{align*}
for
\begin{align*}
v_1(t)&:=\left[t(t+\tfrac{1}{2})(t+1)-(t+1)^3\right]e^{-2\rho t(t+1)},\\
v_2(t)&:=(t+1)^3\left[e^{-2\rho t(t+1)}-e^{-2\rho(t+1)^2}\right],\\
v_3(t)&:=\left[(t+1)^3-(t+1)^2\sqrt{t(t+2)}\right]e^{-2\rho(t+1)^2},
\end{align*}
and for each of the functions~$v_1$, $v_2$, $v_3$ apply the trapezoidal quadrature rule from Lemma~\ref{lem:quadrature} (note, however, that $v_1$ is a negative function, monotonous\-ly decreasing to a minimum, and then monotonously increasing with limit equal to~$0$ in infinity).

It follows from
$$
\int v_1(t)\,dt
=\frac{e^{-2\rho t(t+1)}\left[2\sqrt\rho\,(6t+7)-e^{\tfrac{1}{2}\rho(2t+1)^2}\sqrt{2\pi}\,(3+\rho)\,\text{erf}\left(\frac{\sqrt\rho\,(2t+1)}{\sqrt 2}\right)\right]}{32\rho\sqrt\rho}+C
$$
that
\begin{align*}
\int_0^\infty v_1(t)\,dt&=-\frac{14\sqrt\rho+e^{\,\rho/2}\sqrt{2\pi}\,(3+\rho)\,\text{erfc}\left(\frac{\sqrt\rho}{\sqrt 2}\right)}{32\rho^{\tfrac{3}{2}}}\\
&=-\frac{3\sqrt\frac{\pi}{2}}{16\rho\sqrt\rho}-\frac{1}{4\rho}+\mathcal O\left(\frac{1}{\sqrt\rho}\right)\qquad\text{for }\rho\to0.
\end{align*}
On the other hand, the derivative
$$
v_1^\prime(t)=-\frac{1}{2}\left[6t+5-2\rho(t+1)(2t+1)(3t+2)\right]e^{-2\rho t(t+1)}
$$
disappears for~$t=t_1$ such that
$$
\rho=\frac{6t_1+5}{2(t_1+1)(2t_1+1)(3t_1+2)}=\frac{1}{2t_1^2}+\mathcal O\left(\frac{1}{t_1^3}\right)\qquad\text{for }t_1\to\infty,
$$
i.e.,
$$
t_1=\frac{1}{\sqrt{2\rho}}+\mathcal O(1)\qquad\text{for }\rho\to0.
$$
Thus,
$$
\min_{t\in\mathbb R_+}v_1(t)=v_1(t_1)=-\frac{3}{4e\rho}+\mathcal O\left(\frac{1}{\sqrt\rho}\right).
$$
Further,
\begin{align*}
\int&v_2(t)\,dt=\frac{1}{32\rho^2}\biggl[4\left[1+2\rho(t+1)^2\right]e^{-2\rho(t+1)^2}\\
&-2\left[2+\rho(4t^2+10t+7)\right]e^{-2\rho t(t+1)}+(3+\rho)\sqrt{2\pi\rho}\,e^{\frac{\rho}{2}}\,\text{erf}\left(\sqrt\frac{\rho}{2}\,(2t+1)\right)\biggr]
\end{align*}
and
\begin{align*}
\int_0^\infty v_2(t)\,dt&=\frac{-4(1+2\rho)e^{-2\rho}+2(2+7\rho)+(3+\rho)\sqrt{2\pi\rho}\,e^{\frac{\rho}{2}}\,\text{erfc}\left(\sqrt\frac{\rho}{2}\right)}{32\rho^2}\\
&=\frac{3\sqrt\frac{\pi}{2}}{16\rho\sqrt\rho}+\frac{1}{4\rho}+\mathcal O\left(\frac{1}{\sqrt\rho}\right).
\end{align*}
The derivative
$$
v_2^\prime(t)=(t+1)^2\left[\left[3-2\rho(t+1)(2t+1)\right]e^{2\rho(t+1)}-3+4\rho(t+1)^2\right]e^{-2\rho(t+1)^2}
$$
disappears when
$$
\left[3-2R(2t+1)\right]e^{2R}-3+4R(t+1)=0,
$$
where $R=\rho(t+1)$, i.e., for
$$
t=\frac{(3-2R)\,e^{2R}-3+4R}{4R\,(e^{2R}-1)}=\frac{1}{R}-\frac{3}{4}+\mathcal O(R)\qquad\text{for }R\to0.
$$
The last equation is satisfied for
$$
t=t_2=\frac{1}{\sqrt\rho}-\frac{7}{8}+\mathcal O(\sqrt\rho)\qquad\text{for }\rho\to0.
$$
Consequently,
$$
\max_{t\in\mathbb R_+}v_2(t)=v_2(t_2)=\frac{2}{e^2\rho}+\mathcal O\left(\frac{1}{\sqrt\rho}\right)\qquad\text{for }\rho\to0.
$$

In order to estimate~$v_3$ note that
$$
(t+1)^3-(t+1)^2\sqrt{t(t+2)}=\frac{t+1}{2}+\mathcal O\left(\frac{1}{t}\right)\qquad\text{for }t\to\infty.
$$
Hence, there exists a constant~$c$ such that
$$
\int_0^\infty v_3(t)\,dt=\underbrace{\int_0^1 v_3(t)\,dt}_{=:I_1}+\underbrace{\int_1^\infty\frac{t+1}{2}\,e^{-2\rho(t+1)^2}dt}_{=:I_2}+R
$$
with
$$
|R|\leq c\int_1^\infty\frac{e^{-2\rho(t+1)^2}}{t}\,dt=:cI_3.
$$
Since
$$
I_1=\int_0^1 v_3(t)\,dt\leq1
$$
independently of~$\rho$,
$$
I_2=\int_1^\infty\frac{t+1}{2}\,e^{-2\rho(t+1)^2}\,dt=\frac{e^{-8\rho}}{8\rho}=\frac{1}{8\rho}+\mathcal O(1)\qquad\text{for }\rho\to0,
$$
and
$$
I_3=\int_1^\infty\frac{e^{-2\rho(t+1)^2}}{t}\,dt\leq\int_0^\infty e^{-\rho t^2}dt=\frac{\sqrt\pi}{2\sqrt\rho},
$$
we have
$$
\left|\int_0^\infty v_3(t)\,dt-\frac{1}{8\rho}\right|\leq\mathcal O\left(\frac{1}{\sqrt\rho}\right)\qquad\text{for }\rho\to0.
$$
Further, the derivative
$$
v_3^\prime(t)=(t+1)\left(t+1-\sqrt{t(t+2)}\right)\left[2-\frac{t+1}{\sqrt{t(t+2)}}-4(t+1)^2\rho\right]e^{-2\rho(t+1)^2}
$$
changes its sign from positive to negative in $t=t_3$ such that
$$
\rho=\frac{2\sqrt{t_3(t_3+2)}-t_3-1}{4(t_3+1)^2\sqrt{t_3(t_3+2)}}.
$$
Since
$$
\frac{2\sqrt{t_3(t_3+2)}-t_3-1}{4(t_3+1)^2\sqrt{t_3(t_3+2)}}=\frac{1}{4t_3^2}+\mathcal O\left(\frac{1}{t_3^3}\right)\qquad\text{for }t_3\to\infty,
$$
we have
$$
t_3=\frac{1}{2\sqrt\rho}+\mathcal O(1)\qquad\text{for }\rho\to0
$$
and
$$
\max_{t\in\mathbb R_+}v_3(t)=v_3(t_3)=\frac{1}{4\sqrt{e\rho}}+\mathcal O(\sqrt\rho)\qquad\text{for }\rho\to0.
$$

The desired estimation follows from
\begin{equation*}\begin{split}
&\left|\sum_{l=1}\left[v_1(l)+v_2(l)+v_3(l)\right]-\int_0^\infty\left[v_1(t)+v_2(t)+v_3(t)\right]dt\right|\\
   &\qquad\leq-\min_{t\in[0,\infty)}v_1(t)+\max_{t\in[0,\infty)}v_2(t)+\max_{t\in[0,\infty)}v_3(t).
\end{split}\end{equation*}
\hfill$\Box$

\begin{bfseries}Remark. \end{bfseries}With more sophisticated methods we probably would be able to find a better bound for the error made by replacing $\sum_{l=1}^\infty\left[v_1(l)+v_2(l)\right]$ by $\int_0^\infty\left[v_1(t)+v_2(t)\right]dt$. However, numerical experiments have shown that it would be of order~$\frac{1}{\rho}$, anyway. Consequently, we are not able to compute the exact value of~$U(\Phi_\rho^G)$, as it was the case for Gauss-Weierstrass kernel in~\cite{LFP03}. On the other hand, sole boundedness of the uncertainty product distinguishes Gauss--Weierstrass wavelet from other ones, see~\cite[Theorem~3.4]{IIN16USW}.

\end{document}